\documentclass[11pt]{article}
\usepackage[table]{xcolor}
\usepackage{graphicx}
\usepackage{amsmath,amsfonts,amssymb}
\usepackage{graphicx}
\usepackage{multirow}
\usepackage{tikz,pgf}
\usepackage{url}
\usepackage{amsmath}
\usepackage{graphicx}
\usepackage{epsfig}
\usepackage{epstopdf}
\usepackage{color}
\usepackage{enumitem}
\def\disp{\displaystyle}

\def\tto{\;{\lower 1pt \hbox{$\rightarrow$}}\kern -10pt
\hbox{\raise 2pt \hbox{$\rightarrow$}}\;}

\def\ra{\rangle}
\def\la{\langle}
\def\ve{\varepsilon}
\def\epsilon{\varepsilon}
\def\B{\Bbb B}
\def\h{\hfill\Box}
\def\R{\Bbb R}

\def\N{\Bbb N}
\def\ox{\bar{x}}

\def\oy{\bar{y}}

\def\h{\hfill\square}
\def\dn{\downarrow}
\def\ph{\varphi}

\def\lm{\lambda}

\def\ph{\varphi}

\def\lm{\lambda}

\setlist[enumerate,1]{itemsep=0.0ex,parsep=0.5ex,label={\rm(\alph*)},leftmargin=*, align=left}
\newcounter{lk}

\usepackage[colorlinks=true]{hyperref}

\begin{document}
\begin{center}
{\sc\bf The Fermat-Torricelli Problem and Weiszfeld's Algorithm in the Light of Convex Analysis}\\[1ex]
{\sc Boris S. Mordukhovich}\footnote{Department of Mathematics, Wayne State University, Detroit, Michigan 48202, USA (boris@math.wayne.edu). Research
of this author was partly  supported by the USA National Science Foundation under grants DMS-1007132 and DMS-1512846, by the USA Air Force Office of Scientific Research grant \#15RT0462, and by the Australian Research Council under Discovery Project DP-190100555.},
{\sc Nguyen Mau Nam}\footnote{Fariborz Maseeh Department of Mathematics and Statistics, Portland State University, Portland, OR
97207, USA (mnn3@pdx.edu). Research of this author was partly supported by the USA National
Science Foundation under grant DMS-1716057.}\\[2ex]
{\bf Dedicated to Christiane Tammer with great respect}
\end{center}
\small{\bf Abstract.} In the early 17th century, Pierre de Fermat proposed the following
problem: given three points in the plane, find a point such
that the sum of its Euclidean distances to the three given points is
minimal. This problem was solved by Evangelista Torricelli and was
named the {\em Fermat-Torricelli problem}. A more general version of the Fermat-Torricelli problem asks for a point that minimizes the sum of the distances to a finite number of given points in $\Bbb R^n$. This is one of the main problems in location science. In this paper we revisit the Fermat-Torricelli problem from both theoretical and numerical viewpoints using some ingredients of convex analysis and optimization.\\[1ex]
{\bf Key words.}  Fermat-Torricelli problem; Weiszfeld's algorithm; Convex analysis.\\[1ex]
\noindent {\bf AMS subject classifications.} 49J52, 49J53, 90C31

\newtheorem{Theorem}{Theorem}[section]
\newtheorem{Proposition}[Theorem]{Proposition}
\newtheorem{Remark}[Theorem]{Remark}
\newtheorem{Lemma}[Theorem]{Lemma}
\newtheorem{Corollary}[Theorem]{Corollary}
\newtheorem{Definition}[Theorem]{Definition}
\newtheorem{Example}[Theorem]{Example}
\renewcommand{\theequation}{\thesection.\arabic{equation}}
\normalsize

\section{Introduction}
\setcounter{equation}{0}

\label{intro}

The Fermat-Torricelli problem asks for a point that minimizes the sum of the distances to three given points in the plane. This problem was proposed by Fermat and solved by Torricelli. Torricelli's solution states as follows: if one of the angles of the triangle formed by the three given points is greater than or equal to $120^\circ$, the corresponding vertex is the solution of the problem. Otherwise, the solution is a unique point inside of the triangle formed by the three points such that each side is seen at an angle of $120^\circ$.  The first numerical algorithm for solving the general Fermat-Torricelli problem was introduced by Weiszfeld in 1937 \cite{w}. The assumptions that guarantee the convergence along with the proof were given by Kuhn in 1972. Kuhn also pointed out an example in which  Weiszfeld's algorithm fails to converge; see \cite{k}. The Fermat-Torricelli problem has attracted great attention from many researchers not only because of its mathematical beauty, but also because of its important applications to the field of facility location. Many generalized versions of the Fermat-Torricelli and several new algorithms have been introduced to deal with generalized Fermat-Torricelli problems as well as to improve  Weiszfeld's algorithm; see, e.g., \cite{bs,b1,n2,bn,mns,mv,ul,vz}. The problem has also been revisited several times from different viewpoints as in \cite{ck,d,p,wf} and the references therein. Among many researchers working on the field of facility location, Christiane Tammer is one of the pioneers with important contributions to the theory, algorithms, and software for solving facility location problems; see \cite{t1,t2,t3,t4} and their bibliographies.

The main goal of this paper is to provide easy access to the problem from both theoretical and numerical aspects by using some tools of convex analysis. These tools are presented in the paper with simple proofs that are understandable for students with basic background in introduction to elementary analysis.

The paper is organized as follows. In Section~2 we establish the existence and uniqueness of the optimal solution for the Fermat-Torricelli problem generated by a finite number of points. We also present proofs of some properties of the optimal solution and derive its explicit representations for the case of three points using the \emph{convex subdifferential}. Various advantages of employing convex analysis when solving the Fermat-Torricelli problem has been revealed in \cite{gn,n2,bn,r1}. Section~3 is devoted to revisiting Kuhn's proof of the convergence of Weiszfeld's algorithm. In this section we mainly follow the line of proving the convergence given by Kuhn \cite{k}, while with involving new ingredients from convex analysis to replace some technical tools in order to make the proof more clear.

Throughout the paper the symbol $\B$ denotes the closed unit ball of $\R^n$, and $\B(\ox; r)$ stands for the closed ball centered at $\ox$ with radius $r$.

\section{Elements of Convex Analysis and Properties of Solutions}\label{sec:1}

In this section we review several important concepts of convex analysis to study the classical Fermat-Torricelli problem as well as its general version mentioned in Section~\ref{intro}. Then we present elementary proofs for some properties of optimal solutions of the problem. More details of convex analysis can be found in the fundamental monograph \cite{r}.

Let $\|\cdot\|$ be the Euclidean norm in $\Bbb R^n$. Given a finite number of distinct points $a_i$ for $i=1,\ldots,m$ in $\Bbb R^n$, define the function
\begin{equation}\label{cost}
\ph(x):=\sum_{i=1}^m\|x-a_i\|.
\end{equation}
The mathematical model of the Fermat-Torricelli problem is as follows:
\begin{equation}\label{ft}
\mbox{\rm minimize }\ph(x)\;\mbox{\rm subject to }\;x\in\Bbb R^n.
\end{equation}
The \emph{weighted version} of this problem can be formulated and treated by a similar way.

Let $f\colon\R^n\to\R$ be a real-valued function. The {\em epigraph} of $f$ is a subset
of $\R^n\times\R$ defined by
\begin{equation*}
\mbox{epi }f:=\big\{(x,\alpha)\in\R^{n+1}\;\big|\;x\in\R^n\;\mbox{ and }\;\alpha\ge f(x)\big\}.
\end{equation*}
The function $f$ is called {\em convex} if we have the inequality
\begin{equation*}
f\big(\lm x+(1-\lm)y\big)\le\lm f(x)+(1-\lm)f(y)\;\mbox{ for all }\;x, y\in \R^n\;\mbox{ and }\;\lm\in(0,1).
\end{equation*}
If this equality becomes strict for $x\ne y$, we say that $f$ is \emph{strictly convex}. It is easy to check by the definitions that $f$ is a convex function on $\Bbb R^n$ if and only if its epigraph is a convex set in $\Bbb R^{n+1}$. We obviously have that the function $\ph$ from (\ref{cost}) is convex on $\R^n$.

\begin{Proposition} Let $f\colon\R^n\to \R$ be a convex function. Then $f$ has a local minimum at $\bar x$ if and only if $f$ has an absolute/global minimum at this point.
\end{Proposition}
{\bf Proof:} We only need to prove the ``only if" implication since the converse is trivial. Suppose that $f$ has a local minimum at $\bar x$. Then there exists a number $\delta>0$ with
\begin{equation*}
f(u)\ge f(\bar x)\;\mbox{ for all }\;u\in\B(\bar x;\delta).
\end{equation*}
For any $x\in\R$ it follows that $x_k:=\Big(1-\dfrac{1}{k}\Big)\bar x+\dfrac{1}{k}x\to\bar x$ as $k\to\infty$. Thus we get $x_k\in\B(\bar x; \delta)$ when $k$ is sufficiently large. It follows furthermore that
\begin{equation*}
f(\bar x)\le f(x_k)\le\Big(1-\dfrac{1}{k}\Big)f(\bar x)+\dfrac{1}{k}f(x),
\end{equation*}
which implies in turn that
\begin{equation*}
\dfrac{1}{k}f(\bar x)\le\dfrac{1}{k} f(x),
\end{equation*}
and hence $f(\bar x)\le f(x)$. It shows that $f$ achieves an absolute minimum at $\bar x$. $\h$

\begin{Proposition} The solution set of the Fermat-Torricelli problem {\rm (\ref{ft})} is nonempty.
\end{Proposition}
{\bf Proof: } Let $m:=\inf\{\ph(x)\;|\;x\in\R^n\}\ge 0$, and let $(x_k)$ be a sequence satisfying
\begin{equation*}
\lim_{k\to\infty}\ph(x_k)=m.
\end{equation*}
It follows from the definition that there exists $k_0\in\Bbb N$ with
\begin{equation*}
\|x_k-a_1\|\le\ph(x_k)\le m+1\;\mbox{\rm for all }\;k\ge k_0,
\end{equation*}
which implies that $\|x_k\|\le m+1+\|a_1\|$ for such $k$. Thus $\{x_k\}$ is a bounded sequence, and so it has a subsequence $\{x_{k_\ell}\}$ that converges to $\ox\in\R^n$ as $\ell\to\infty$. Since $\ph$ is obviously continuous, we get
\begin{equation*}
\ph(\ox)=\lim_{\ell\to\infty}\ph(x_{k_\ell})=m.
\end{equation*}
This shows that $\ox$ is an optimal solution of the problem. $\h$\vspace*{0.05in}

For two different points $a,b\in\R^n$, the line containing $a$ and $b$ is given by
\begin{equation*}
\mathcal{L}(a,b):=\big\{ta+(1-t)b\;\big|\;t\in\R\big\}.
\end{equation*}

\begin{Proposition} Suppose that the points $a_i$ for $i=1,\ldots,m$ do not lie on the same line (i.e., not collinear). Then the function $\ph$ defined by {\rm (\ref{cost})} is strictly convex, and the Fermat-Torricelli problem {\rm(\ref{ft})} has a unique solution.
\end{Proposition}
{\bf Proof: }Define $\ph_i(x):=\|x-a_i\|$ for $i=1,\ldots,m$. Then $\ph=\sum_{i=1}^m\ph_i$. For any $x,y\in\R^n$ and $\lambda\in(0,1)$ we have the inequality
\begin{equation*}
\ph_i\big(\lambda x+(1-\lambda)y\big)\le\lambda\ph_i(x)+(1-\lambda)\ph_i(y)\;\mbox{\rm for }\;i=1,\ldots,m.
\end{equation*}
This readily implies that
\begin{equation}\label{convexity}
\ph\big(\lambda x+(1-\lambda)y\big)\le\lambda\ph(x)+(1-\lambda)\ph(y).
\end{equation}
On the contrary, suppose that $\ph$ is not strictly convex. It means that there exist $\ox,\oy\in\R^n$ with $\ox\ne\oy$ and $\lambda\in(0,1)$ for which (\ref{convexity}) holds as equality. Then
\begin{equation*}
\ph_i\big(\lambda\ox+(1-\lambda)\oy\big)=\lambda\ph_i(\ox)+(1-\lambda)\ph_i(\oy)\;\mbox{\rm for }\;i=1,\ldots,m,
\end{equation*}
which can be rewritten as
\begin{equation*}
\|\lambda (\ox-a_i)+(1-\lambda)(\oy-a_i)\|=\|\lambda(\ox-a_i)\|+\|(1-\lambda)(\oy-a_i)\|\;\mbox{\rm for }\;i=1,\ldots,m.
\end{equation*}
If $\ox\ne a_i$ and $\oy\ne a_i$, then there exists $t_i>0$ such that
\begin{equation*}
t_i\lambda(\ox-a_i)=(1-\lambda)(\oy-a_i).
\end{equation*}
This tells us that $\ox-a_i=\gamma_i(\oy-a_i)$, where $\gamma_i:=\dfrac{1-\lambda}{t_i\lambda}$. Since $\ox\ne\oy$, we obtain $\gamma_i\ne 1$ and
\begin{equation*}
a_i=\dfrac{1}{1-\gamma_i}\ox-\dfrac{\gamma_i}{1-\gamma_i}\oy\in\mathcal{L}(\ox,\oy).
\end{equation*}
In the case where either $\ox=a_i$ or $\oy=a_i$, it is obvious that $a_i\in\mathcal{L}(\ox,\oy)$. Thus we arrive at $a_i\in\mathcal{L}(\ox,\oy)$ for $i=1,\ldots,m$, which is a contradiction that completes the proof. $\h$\vspace*{0.05in}

Recall that a vector $v\in\R^n$ is a {\em subgradient} of a convex function $f\colon\R^n\to\R$ at the given point $\ox\in\R^n$ if it satisfies the inequality
\begin{equation}\label{sub}
f(x)\ge f(\ox)+\la v,x-\ox\ra\;\mbox{ for all }\;x\in\R^n,
\end{equation}
where $\la\cdot,\cdot\ra$ stands for the usual scalar product in $\R^n$. The set of all subgradients of $f$ at $\ox$ is called the \emph{subdifferential} of this function at $\ox$ and is denoted by $\partial f(\ox)$.

Directly from the definition, we deduce the following {\em subdifferential Fermat rule}:
\begin{equation}\label{fermat}
\mbox{\rm $f$ achieves an absolute minimum at $\ox$ if and only if $0\in\partial f(\ox)$.}
\end{equation}
The proposition below shows that the subdifferential of a convex function at a given point reduces to the gradient at that point when the function is differentiable.

\begin{Proposition}\label{dd1} Suppose that $f\colon\R^n\to\R$ is convex and (Fr\'echet) differentiable at $\ox$. Then
\begin{equation}\label{d1}
\la\nabla f(\ox),x-\ox\ra\le f(x)-f(\ox)\;\mbox{ for all }\;x\in\R^n.
\end{equation}
Furthermore, we have $\partial f(\ox)=\{\nabla f(\ox)\}$.
\end{Proposition}
{\bf Proof: }Since $f$ is differentiable at $\ox$, it follows from the definition that for any $\epsilon>0$ there exists $\delta>0$ such that
\begin{equation*}
-\epsilon\|x-\ox\|\le f(x)-f(\ox)-\la\nabla f(\ox),x-\ox\ra\le\epsilon\|x-\ox\|\;\mbox{\rm whenever }\;\|x-\ox\|<\delta.
\end{equation*}
Consider further the function
\begin{equation*}
\psi(x):=f(x)-f(\ox)-\la\nabla f(\ox),x-\ox\ra+\epsilon\|x-\ox\|,
\end{equation*}
which satisfies $\psi(x)\ge\psi(\ox)=0$ for all $x\in\B(\ox;\delta)$. By the convexity of $\psi$ we have $\psi(x)\ge\psi(\ox)$ whenever $x\in\R^n$. Hence
\begin{equation*}
\la\nabla f(\ox),x-\ox\ra\le f(x)-f(\ox)+\epsilon\|x-\ox\|\;\mbox{ for all }\;x\in\R^n.
\end{equation*}
Letting now $\epsilon\dn 0$ justifies the estimate in (\ref{d1}).

Equality (\ref{d1}) yields $\nabla f(\ox)\in\partial f(\ox)$. Taking any $v\in\partial f(\ox)$, we get by definition that
\begin{equation*}
\la v,x-\ox\ra\le f(x)-f(\ox)\;\mbox{\rm for all }\;x\in\R^n.
\end{equation*}
The differentiability of $f$ at $\ox$ also implies that for any $\epsilon>0$ there exists $\delta>0$ such that
\begin{equation*}
\la v-\nabla f(\ox),x-\ox\ra\le\epsilon\|x-\ox\|\;\mbox{ whenever }\;\|x-\ox\|<\delta.
\end{equation*}
Thus we have $\|v-\nabla f(\ox)\|\le\epsilon$, which shows that $v=\nabla f(\ox)$ since $\epsilon>0$ was chosen arbitrarily. It verifies that $\partial f(\ox)=\{\nabla f(\ox)\}$ and thus completes the proof of the proposition. $\h$\vspace*{0.05in}

The subdifferential formula for the norm function derived in the next example plays a crucial role in our subsequent analysis to solve the Fermat-Torricelli problem.

\begin{Example} {\rm Let $p(x):=\|x\|$, the Euclidean norm function on $\R^n$. Then we have
\begin{equation*}
\partial p(x)=\begin{cases}
\B &\text{if }\;x=0,\\
\Big\{\dfrac{x}{\|x\|}\Big\}&\text{otherwise}.
\end{cases}
\end{equation*}
Since the function $p(\cdot)$ is differentiable with $\nabla p(x)=x/\|x\|$ for $x\ne 0$, it suffices to verify the claimed formula for $x=0$. It follows from \eqref{sub} that $v\in\partial p(0)$ if and only if
\begin{equation*}
\la v,x\ra=\la v,x-0\ra\le p(x)-p(0)=\|x\|\;\mbox{\rm for all }\;x\in\R^n.
\end{equation*}
Letting $x:=v$, we get $\la v,v\ra\le\|v\|$, which implies that $\|v\|\le 1$, i.e., $v\in\B$. Conversely, for $v\in\B$ the Cauchy-Schwarz inequality tells us that
\begin{equation*}
\la v,x-0\ra=\la v,x\ra\le\|v\|\cdot\|x\|\le\|x\|=p(x)-p(0)\;\mbox{\rm whenever }\;x\in\R^n,
\end{equation*}
and hence $v\in\partial p(0)$. Thus we arrive at $\partial p(0)=\B$.}
\end{Example}

Solving the Fermat-Torricelli problem involves the usage of the following subdifferential rule for sums of two convex functions one of which is differentiable while the other may be not.

\begin{Proposition}\label{sr} Let $f_i\colon\R^n\to\R$, $i=1,2$, be two convex functions such that $f_2$ is differentiable at $\ox$. Then we have the equality
\begin{equation}\label{srl}
\partial(f_1+f_2)(\ox)=\partial f_1(\ox)+\nabla f_2(\ox).
\end{equation}
\end{Proposition}
{\bf Proof: }Fix any $v\in\partial(f_1+f_2)(\ox)$ and get for each $x\in\R^n$ that
\begin{equation*}
\la v,x-\ox\ra\le f_1(x)-f_1(\ox)+f_2(x)-f_2(\ox)=f_1(x)-f_1(\ox)+\la\nabla f_2(\ox),x-\ox\ra+o(\|x-\ox\|).
\end{equation*}
For any $\epsilon>0$ there exists $\delta>0$ such that
\begin{equation*}
0\le\la\nabla f_2(\ox)-v,x-\ox\ra+f_1(x)-f_1(\ox)+\epsilon\|x-\ox\|\;\mbox{\rm whenever }\;x\in\B(\ox;\delta).
\end{equation*}
The convexity of $f_1$ ensures that the latter holds for all $x\in\R^n$. Letting $\epsilon\dn0$ yields
\begin{equation*}
0\le\la\nabla f_2(\ox)-v,x-\ox\ra+f_1(x)-f_1(\ox)\;\mbox{\rm whenever }\;x\in\R^n.
\end{equation*}
By \eqref{sub} it tells us that $v-\nabla f_2(\ox)\in\partial f_1(\ox)$ and hence $v\in\partial f_1(\ox)+\nabla f_2(\ox)$, which justifies the inclusion ``$\subset$" in (\ref{srl}). Representing
$f_1=f_1+f_2+(-f_2)$ and applying the obtained inclusion, we have the relationships
\begin{equation*}
\partial f_1(\ox)\subset\partial(f_1+f_2)(\ox)+\nabla(-f_2)(\ox)=\partial(f_1+f_2)(\ox)-\nabla f_2(\ox),
\end{equation*}
which verify the opposite inclusion ``$\supset$" in (\ref{srl}) and thus complete the proof. $\h$\vspace*{0.05in}

Let us now use subgradients of the norm function to obtain the classical solution of the Fermat-Torricelli problem. Given two nonzero vectors $u,v\in\R^n$, denote
\begin{equation*}
\cos(u,v):=\dfrac{\la u,v\ra}{\|u\|\cdot\|v\|}.
\end{equation*}
Fix $\ox\ne a_i$ and define the unit vectors
\begin{equation*}\label{v}
v_i:=\dfrac{\ox-a_i}{\|\ox-a_i\|},\;i=1,2,3.
\end{equation*}
Each $v_i$ is the unit vector pointing in the direction from the vertex $a_i$ to $\ox$. Observe that the Fermat-Torricelli problem formulated above always has a unique solution even if the three given points are on the same line. It is easy to see in the latter case that the middle point is the solution of the problem. The next proposition completely characterizes the solution of the Fermat-Torricelli problem in the general three-point setting of the $n$-dimensional space, not just on the plane as in the original framework.

\begin{Proposition}\label{3point} Let the points $a_1,a_2,a_3\in\R^n$ generate the Fermat-Torricelli problem in $\R^n$. Then we have the following descriptions of the optimal solution to \eqref{ft} with $\ph$ taken from \eqref{cost}:\\
{\bf(i)} In the case where $\ox\notin\{a_1,a_2,a_3\}$, $\ox$ is the solution of the problem if and only if
\begin{equation*}
\cos(v_1,v_2)=\cos(v_2,v_3)=\cos(v_3,v_1)=-1/2.
\end{equation*}
{\bf(ii)} Consider the case where $\ox\in\{a_1,a_2,a_3\}$ and suppose for definiteness that $\ox=a_1$. Then $\ox$ is the solution of the problem if and only if
\begin{equation*}
\cos\la v_2,v_3\ra\le-1/2.
\end{equation*}
\end{Proposition}
{\bf Proof: }In case (i) we have that the function $\ph$ from (\ref{cost}) is differentiable at $\ox$. Since $\ph$ is convex, $\ox$ is the solution to the Fermat-Torricelli problem if and only if
\begin{equation*}
\nabla\ph(\ox)=v_1+v_2+v_3=0.
\end{equation*}
Remembering that $\|v_i\|=1$ for $i=1,2,3$, we get
\begin{align*}
&\la v_1,v_2\ra+\la v_1,v_3\ra=-1\\
&\la v_2,v_1\ra+\la v_2,v_3\ra=-1\\
&\la v_3,v_1\ra+\la v_3,v_2\ra=-1.
\end{align*}
Solving this system of equations yields
\begin{equation*}
\la v_i,v_j\ra =\cos(v_i,v_j)=-1/2\;\mbox{\rm for }\;i\ne j,\;i,j\in\{1,2,3\}.
\end{equation*}
If furthermore $\la v_i,v_j\ra =-1/2$ for $i\ne j$, $i,j\in\{1,2,3\}$, then
\begin{equation*}
\|v_1+v_2+v_3\|^2=\sum_{i=1}^3\|v_i\|+\sum_{i,j=1, i\neq j}^3\la v_i,v_j\ra =0,
\end{equation*}
which gives us $v_1+v_2+v_3=0$ and thus completes the proof in the case.\\[1ex]
In case (ii) we deduce from the subdifferential Fermat rule (\ref{fermat}) and the subdifferential sum rule (\ref{srl}) that $\ox=a_1$ is the solution to the Fermat-Torricelli problem if and only if
\begin{equation*}
0\in\partial\ph(a_1)=\B+v_2+v_3.
\end{equation*}
This is equivalent to $\|v_2+v_3\|^2\le 1$ or, equivalently, to $\|v_2\|^2+\|v_3\|^2+2\la v_2,v_3\ra\le 1$. Since $v_2$ and $v_3$ are unit vectors, we obtain
\begin{equation*}
\la v_2,v_3\ra=\cos(v_2,v_3)\le-1/2
\end{equation*}
and complete the proof of the proposition. $\h$\vspace*{0.05in}

Next we present an example with a figure illustrating the obtained solution in the classical case of three points on the plane.

\begin{Example}{\rm Consider the Fermat-Torricelli problem given by three points $A$, $B$, and $C$ on the plane as shown in the figure. If one of the angles of the triangle $ABC$ is greater than or equal to $120^\circ$, then the corresponding vertex is the solution to the problem by Proposition~\ref{3point}(ii). Let us examine the case where none of the angles of the triangle is greater than or equal to $120^\circ$. Construct two equilateral triangles $ABD$ and $ACE$ and let $S$ be the intersection of $DC$ and $BE$ as in the figure. Two quadrilaterals $ADBC$ and $ABCE$ are convex, and hence $S$ lies inside the triangle $ABC$.  It is clear that two triangles $DAC$ and $BAE$ are congruent. A rotation of $60^\circ$ about $A$ maps the triangle $DAC$ to the triangle $BAE$. The rotation maps $CD$ to $BE$, so $\angle DSB=60^\circ$. Let $T$ be the image of $S$ through this rotation. Then $T$ belongs to $BE$. It follows that $\angle AST=\angle ASE=60^\circ$. Moreover, $\angle DSA=60^\circ$, and hence $\angle BSA=120^\circ$. It is now clear that $\angle ASC=120^\circ$ and $\angle BSC=120^\circ$. Proposition~\ref{3point}(i) tells us that the point $S$ is the solution to this classical Fermat-Torricelli problem.}

\begin{center}
\includegraphics[scale=0.70]{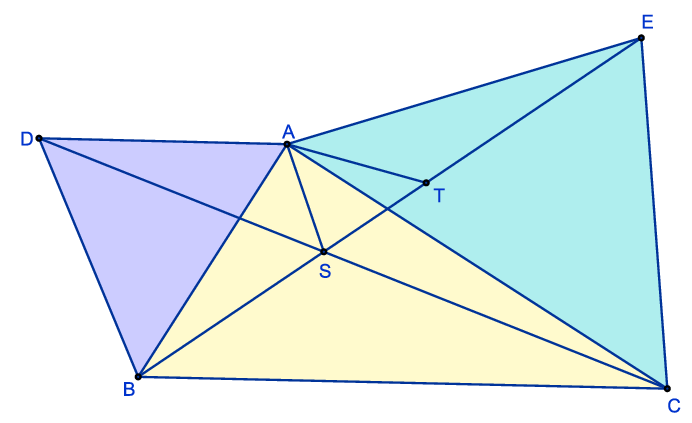}
\end{center}
\end{Example}

\section{Weiszfeld's Algorithm}

In this section we revisit Kuhn's proof \cite{k} of the convergence of Weiszfeld's algorithm \cite{w} for solving the Fermat-Torricelli problem (\ref{ft}). With some additional ingredients of convex analysis we are able to provide a more clear picture of Kuhn's proof. Throughout this section, assume that the points $a_i$ for $i=1,\ldots,m$ are not collinear.

The gradient of the function $\ph$ from (\ref{cost}) is computed by
\begin{equation*}
\nabla\ph(x)=\sum_{i=1}^m\dfrac{x-a_i}{\|x-a_i\|},\; x\notin\big\{a_1,a_2,\ldots,a_m\big\}.
\end{equation*}
Solving the gradient equation $\nabla\ph(x)=0$ gives us the formula
\begin{equation}\label{grad}
x=\dfrac{\disp\sum_{i=1}^m\dfrac{a_i}{\|x-a_i\|}}{\disp\sum_{i=1}^m\dfrac{1}{\|x-a_i\|}}:=F(x).
\end{equation}
To keep the continuity, define $F(x):=x$ for $x\in\{a_1,a_2,\ldots,a_m\}$.

Weiszfeld introduced the following algorithm: choose a starting point $x_0\in\R^n$ and define
\begin{equation*}
x_{k}=F(x_{k-1})\;\mbox{\rm for }\;k\in\N.
\end{equation*}
He also claimed that if $x_0\notin\{a_1,a_2,\ldots,a_m\}$ where $a_i$ for $i=1,\ldots,m$ are not collinear, then $\{x_k\}$ converges to the unique optimal solution of the problem. A correct statement and the proof of the convergence were given by Kuhn in \cite{k}.\vspace*{0.05in}

The next proposition guarantees that the cost function value in \eqref{ft} decreases after each iteration of the Weiszfeld algorithm.

\begin{Proposition}\label{dp} If $F(x)\ne x$, then $\ph(F(x))<\ph(x)$.
\end{Proposition}
{\bf Proof: }It is clear that $x$ is not a vertex, since otherwise we get $F(x)=x$. Moreover, the point $F(x)$ is the unique minimizer of the following strictly convex function:
\begin{equation*}
g(z):=\sum_{i=1}^m\dfrac{\|z-a_i\|^2}{\|x-a_i\|}.
\end{equation*}
Indeed, we can easily check that $F(x)$ is the unique solution of the equation $\nabla g(z)=0$. Using $F(x)\ne x$ tells us that $g(F(x))<g(x)=\ph(x)$. Furthermore, it clearly follows from the constructions above that
\begin{align*}
g\big(F(x)\big)&=\sum_{i=1}^m\dfrac{\|F(x)-a_i\|^2}{\|x-a_i\|}\\
&=\sum_{i=1}^m\dfrac{(\|x-a_i\|+\|F(x)-a_i\|-\|x-a_i\|)^2}{\|x-a_i\|}\\
&=\ph(x)+2\big(\ph(F(x))-\ph(x)\big)+\sum_{i=1}^m\dfrac{(\|F(x)-a_i\|-\|x-a_i\|)^2}{\|x-a_i\|},
\end{align*}
which verifies the strict inequality
\begin{equation*}
2\ph\big(F(x)\big)+\sum_{i=1}^m\dfrac{(\|F(x)-a_i\|-\|x-a_i\|)^2}{\|x-a_i\|}<2\ph(x)
\end{equation*}
and hence yields the claimed decreasing property $\ph(F(x))<\ph(x)$. $\h$\vspace*{0.05in}

Now we investigate behavior of the algorithm mapping $F$ near a vertex and deal with the case where a vertex is not the solution of the Fermat-Torricelli problem \eqref{ft}. Let us first present a necessary and sufficient condition for a vertex to be the optimal solution of the problem. Define
\begin{equation*}
R_j:=\sum_{i=1,i\ne j}^m\dfrac{a_i-a_j}{\|a_i-a_j\|},\quad j=1,\ldots,m,
\end{equation*}

\begin{Proposition}\label{lm2} The vertex $a_j$ is the optimal solution to \eqref{ft} if and only if $\|R_j\|\le 1$.
\end{Proposition}
{\bf Proof: }Employing the subdifferential Fermat rule (\ref{fermat}) and the subdifferential sum rule from Proposition~\ref{sr} ensures that the vertex $a_j$ is the optimal solution of the problem if and only if
\begin{equation*}
0\in\partial\ph(a_j)=-R_j+\B,
\end{equation*}
which can be equivalently rewritten as $\|R_j\|\le 1$. $\h$\vspace*{0.05in}

The obtained result allows us to significantly simplify the proof of the next proposition taken from \cite[Subsection~3.2]{k}.

\begin{Proposition}\label{k} Suppose that the vertex $a_j$ is not the optimal solution to \eqref{ft}. Then there is a number $\delta>0$ such that the condition $0<\|x-a_j\|\le\delta$ yields the existence of a positive integer $q$ for which we have the estimates
\begin{equation}\label{dd}
\|F^q(x)-a_j\|>\delta\;\mbox{\rm and }\;\|F^{q-1}(x)-a_j\|\le\delta
\end{equation}
with using the notation
\begin{equation}\label{not}
F^q(x):=F\big(F^{q-1}(x)\big)\;\mbox{ whenever }\;q=1,2,\ldots\;\mbox{ and }\;F^0(x):=x.
\end{equation}
\end{Proposition}
{\bf Proof: }If $x$ is not a vertex, then we get from \eqref{grad} that
\begin{align*}
F(x)=\dfrac{\disp\sum_{i=1}^m\dfrac{a_i}{\|x-a_i\|}}{\disp\sum_{i=1}^m\dfrac{1}{\|x-a_i\|}}.
\end{align*}
which implies in turn that
\begin{equation*}
F(x)-a_j=\dfrac{\disp\sum_{i=1,i\ne j}^m\dfrac{a_i-a_j}{\|x-a_i\|}}{\disp\sum_{i=1}^m\dfrac{1}{\|x-a_i\|}}.
\end{equation*}
Taking now the limit as $x\to a_j$ leads us to
\begin{equation*}
\lim_{x\to a_j}\dfrac{F(x)-a_j}{\|x-a_j\|}=\lim_{x\to a_j}\dfrac{\disp\sum_{i=1,i\ne j}^m\dfrac{a_i-a_j}{\|x-a_i\|}}{1+\disp\sum_{i=1,i\ne j}^m\dfrac{\|x-a_j\|}{\|x-a_i\|}}=R_j.
\end{equation*}
This implies by Proposition~\ref{lm2} that
\begin{equation}\label{R}
\lim_{x\to a_j}\dfrac{\|F(x)-a_j\|}{\|x-a_j\|}=\|R_j\|>1
\end{equation}
and thus allows us to find positive numbers $\epsilon$ and $\delta$ with
\begin{equation}\label{ep}
\dfrac{\|F(x)-a_j\|}{\|x-a_j\|}\geq (1+\epsilon)\;\mbox{\rm whenever }0<\|x-a_j\|\leq\delta.
\end{equation}
Remembering the notation in \eqref{not}, if
$$
0<\|F^{p-1}(x)-a_j\|\le\delta\;\mbox{ for all }\;p=1,\ldots,q,
$$
then by \eqref{ep} we have
$$
\|F^q(x)-a_j\|\ge(1+\ve)\|F^{q-1}(x)-a_j\|\ge\ldots\ge(1+\ve)^q\|x-a_j\|.
$$
Taking into account that $(1+\ve)^q\|x-a_j\|\to\infty$ as $q\to\infty$ verifies the estimates in \eqref{dd} and thus completes the proof of the proposition. $\h$\vspace*{0.05in}

We finally present the following simplified and improved proof (with taking into account Propositions~\ref{lm2} and \eqref{k} above) of Kuhn's convergence result \cite{k} for Weiszfeld's algorithm to solve the Fermat-Torricelli problem.

\begin{Theorem} Let $\{x_k\}$ be the sequence of iterates generated by Weiszfeld's algorithm, and let $x_k\notin\{a_1,a_2,\ldots,a_m\}$ for all $k=0,1,\ldots$. Then $\{x_k\}$ converges to the optimal solution $\ox$ of the Fermat-Torricelli problem \eqref{ft}.
\end{Theorem}
{\bf Proof: }Observe first that if $x_k=x_{k+1}$ for some $k=k_0$, then $x_k$ is a constant sequence for all $k\ge k_0$, which therefore converges to $x_{k_0}$. Since $F(x_{k_0})=x_{k_0}$ and $x_{k_0}$ is not a vertex, the point $x_{k_0}$ is the solution of the problem. Hence we can proceed by assuming that $x_{k+1}\ne x_k$ for every $k$. Proposition~\ref{dp} tells us that the sequence $\{\ph(x_k)\}$ is nonnegative and decreasing, and thus it converges, which means that
\begin{equation}\label{dc}
\lim_{k\to\infty}\big(\ph(x_k)-\ph(x_{k+1})\big)=0.
\end{equation}
It follows from the algorithm that $x_k\in\mbox{\rm co }\{a_1,a_2,\ldots,a_m\}$ (the convex hull) for all $k\geq 1$. Since the latter set is compact in $\R^n$, we have the convergence of some subsequence of $\{x_k\}$. Take a subsequence $\{x_{k_\ell}\}$ of $\{x_k\}$ that converges to a point $\oy$. It suffices to prove that $\oy=\ox$. To proceed, deduce from (\ref{dc}) that
\begin{equation*}
\lim_{\ell\to\infty}\big(\ph(x_{k_\ell})-\ph(F(x_{k_\ell})\big)=0
\end{equation*}
and conclude by the continuity of $\ph$ that $\ph(\oy)=\ph(F(\oy))$. This clearly yields $F(\oy)=\oy$.

If $\oy$ is not a vertex, then it is the solution of the problem, so $\oy=\ox$. Let us consider the case where $\oy$ is a vertex, say $a_1$. Arguing by contradiction, suppose  that $\oy\ne\ox$. Choose $\delta>0$ sufficiently small such that the properties in Proposition~\ref{k} hold, and that the ball $\B(a_1;\delta)$ does not contain $\ox$ and $a_i$ for $i=2,\ldots,m$. Since $x_{k_\ell}\to a_1=\oy$, we assume without loss of generality that the sequence is contained in $\B(a_1;\delta)$.

For $x=x_{k_1}$, choose $q_1$ such that $x_{q_1}\in\B(a_1;\delta)$ and $F(x_{q_1})\notin\B(a_1;\delta)$. Selecting further an index $k_\ell>q_1$ and applying Proposition~\ref{k}, we find $q_2>q_1$ such that $x_{q_2}\in\B(a_1;\delta)$ and $F(x_{q_2})\notin\B(a_1;\delta)$. Repeating this procedure gives us a sequence $\{x_{q_\ell}\}$ with $x_{q_\ell}\in\B(a_1;\delta)$ and $F(x_{q_\ell})$ not belonging to the ball. Extracting yet another subsequence, suppose that $x_{q_\ell}\to\bar z$. It follows from the above that $F(\bar z)=\bar z$. If $\bar z$ is not a vertex, then it must be the solution, which is a contradiction because the solution $\ox$ is not in $\B(a_1;\delta)$. Thus $\bar z$ is a vertex that should be $a_1$, since the other vertices do not belong to the ball as well. It tells us that
\begin{equation*}
\lim_{\ell\to\infty}\dfrac{\|F(x_{q_\ell})-a_1\|}{\|x_{q_\ell}-a_1\|}=\infty,
\end{equation*}
which contradicts Proposition~\ref{lm2} via \eqref{R} and thus completes the proof of the theorem.$\h$

\vskip 6mm
\noindent{\bf Acknowledgments.} \noindent Research of Boris Mordukhovich was partly supported by the USA National Science Foundation under grants DMS-1007132 and DMS-1512846, by the USA Air Force Office of Scientific Research grant \#15RT0462, and by the Australian Research Council under Discovery Project DP-190100555. Research of Nguyen Mau Nam was partly supported by the USA National Science Foundation under grant DMS-1716057.


\begin{thebibliography}{99}

\bibitem{t1} S. Alzorba, C. Gunther, N. Popovici and C. Tammer,  A new algorithm for solving planar multi-objective location problems involving the Manhattan norm. European
J. Oper. Res. 258 ( 2017), 35--46.

\bibitem{bs} A. Beck and S. Sabach, Weiszfeld's method: old and new results. J. Optim. Theory Appl. 164 (2015), 1--40.

\bibitem{b} J. Brimberg. The Fermat Weber location problem revisited. Math. Program. 71 (1995), 71--76.

\bibitem{b1} J. Brimberg, R. Chen and D. Chen, Accelerating convergence in the Fermat-Weber location problem. Oper. Res. Lett. 22 (1998), 151--157.

\bibitem{ck} L. Cooper and I. Katz, The Weber problem revisited. Comput. Math. Appl. 7 (1981), 225--234.

\bibitem{d} Z. Drezner, On the convergence of the generalized Weiszfeld algorithm. Ann. Oper. Res. 167 (2009), 327--336.

\bibitem{e} U. Eckhardt, Weber's problem and Weiszfeld's algorithm in general spaces. Math. Program. 18 (1980) 186--196.

\bibitem{gn} F. Giannessi, Constrained Optimization and Image Space Analysis, Vol. 1. Separation of Sets and Optimality Conditions. Math. Concepts Methods Sci. Engrg. 49, Springer, New York (2005).

\bibitem{t2} C. Gunther and C. Tammer, Relationsship between constrained and unconstrained multi-objective optimization and application in location theory. Math. Methods Oper. Res. 84 (2016), 359--387.

\bibitem{t3} C. Gunther, C. Tammer, M. Hillmann and B. Winkler, Project Facility Location Optimizer. The Martin Luther University of Halle-Wittenberg (2018).

\bibitem{k} H.W. Kuhn, A note on Fermat-Torricelli problem. Math. Program. 4 (1973), 98--107.

\bibitem{n2} B.S. Mordukhovich and N.M. Nam, Applications of variational analysis to a generalized Fermat-Torricelli problem. J. Optim. Theory Appl.
148 (2011), 431--454.

\bibitem{bn} B.S. Mordukhovich and N. M. Nam, An Easy Path to Convex Analysis and Applications. Morgan \& Claypool Publishers, San Rafael, CA (2014).

\bibitem{mns} B.S. Mordukhovich, N.M. Nam and J. Salinas, Applications of variational analysis to a generalized Heron problem. Applic Anal. 91 (2012), 1915--1942.

\bibitem{mv} J.G. Morris and W.A. Verdini, Minisum $\ell_p$ distance location problems solved via a perturbed problem and Weiszfeld's algorithm. Oper. Res. 27 (1979), 1180--1188.

\bibitem{t4} R. Patz, J. Spitzner and C. Tammer, Decision support for location problems in town planning. Inter. Trans. Oper. Res. 9 (2002), 261--278.

\bibitem{p} F. Plastria, The Weiszfeld algorithm: proof, amendments and extensions. H.A. Eiselt and V. Marianov (Eds.). Foundations of Location Analysis. International Series in Operations Research
and Management Science, Vol. 155, pp.\ 357--389. Springer, New York (2011).

\bibitem{r} R.T. Rockafellar, Convex Analysis. Princeton University Press, Princeton, NJ (1970).

\bibitem{r1} A. Ruszczy\'nski, Nonlinear Optimization. Princeton University Press, Princeton, NJ, 2006.

\bibitem{ul} H. $\ddot{\rm U}$ster and R.F. Love, The convergence of the Weiszfeld algorithm. Comput. Math. Appl. 40 (2000), 443--451.

\bibitem{vz} Y. Vardi and C.-H. Zhang,  A modified Weiszfeld algorithm for the Fermat-Weber location problem. Math. Program. 90 (2001), 559--566.

\bibitem{w} E. Weiszfeld, Sur le point pour lequel la somme des distances de $n$ points donn\'es est minimum, T$\hat{\mbox{\rm o}}$hoku Mathematics Journal 43 (1937), 355--386.

\bibitem{wf} E. Weiszfeld and F. Plastria, On the point for which the sum of the distances to n given points is minimum. Ann Oper Res 167 (2009), 7--41.

\end{thebibliography}
\end{document}